\documentclass[12pt,reqno]{amsart}
\usepackage{latexsym, amsfonts, amsmath, amsthm, amssymb,amscd}
\makeatletter
\newfont{\footsc}{cmcsc10 at 8truept}
\newfont{\footbf}{cmbx10 at 8truept}
\newfont{\footrm}{cmr10 at 10truept}
\makeatother 

\newtheorem{thm}{Theorem}[section]
\newtheorem{prop}[thm]{Proposition}

\newtheorem{cor}[thm]{Corollary}
\newtheorem{defn}[thm]{Definition}
\newtheorem{lem}[thm]{Lemma}

\newtheorem{exam}{Example}[section]
\newtheorem{rmk}[thm]{Remark}

\numberwithin{equation}{section}

\def\pf{\noindent {\it Proof.} }

\def\nV{\hbox{${}\kern-2pt\not\kern-2pt {\rm V}$}}
\def\nU{\hbox{${}\kern-2pt\not\kern-2pt {\rm U}$}}
\def\nS{\hbox{${}\kern-2pt\not\kern-2pt {\rm S}$}}
\def\nH{\hbox{${}\kern-2pt\not\kern-2pt {\rm H}$}}
\def\nT{\hbox{${}\kern-2pt\not\kern-2pt {\rm T}$}}
\def\V{{\rm V}}
\def\U{{\rm U}}
\def\H{{\rm H}}
\def\S{\rm S}
\def\T{\rm T}
\def\R{{\mathcal R}}

\def\w{{\rm w}}
\def\c{\textbf{c}}

\def\inv{\mathop {\rm inv}}
\def\maj{\mathop {\rm maj}}
\def\stat{\mathop {\rm stat}}

\newcommand\MAJ{\operatorname{MAJ}}

\title{A classification of mahonian maj-inv statistics}

\author{Anisse Kasraoui}

\begin{document}
\maketitle

\centerline{\small Universit\'e de Lyon;}
 \centerline{\small Universit\'e Lyon 1;}
 \centerline{\small CNRS, UMR5208, Institut Camille Jordan}
 \centerline{INSA de Lyon, F-69621;}
  \centerline{Ecole centrale de Lyon;}
 \centerline{\small 43, boulevard du 11 novembre 1918,}
 \centerline{F-69622 Villeurbanne Cedex, France}

\vspace{0.3cm}
 \centerline{\small\texttt{anisse@math.univ-lyon1.fr}}
\centerline{\small\texttt{tel: +33 (0)4 72 43 11 89 }}
\centerline{\small\texttt{fax: +33 (0)4 72 43 16 87 }}

\vspace{0.3cm}

\begin{abstract}
Two well-known mahonian statistics on words  are the inversion
number and the major index. In 1996, Foata and Zeilberger introduced
generalizations, parametrized by relations, of these statistics. In
this paper, we study the statistics which can be written as a sum of
these generalized statistics. This leads to generalizations of some
classical results. In particular, we characterize all such
statistics which are mahonian.
\end{abstract}

\vspace{0.5cm}

\vspace{0.5cm}

\noindent{\bf MR Subject
Classifications}: Primary 05A05, 05A15; Secondary 05A30.\\

\section{Introduction and main results}

\subsection{Introduction} Let $X$ be a finite alphabet. Without loss of generality
we may assume $X=[r]:=\{1,2,\ldots,r\}$. Two of the most known and
studied statistics on words (and permutations) are probably the
\emph{inversion number} ($\inv$) and the \emph{major index}
($\maj$). They are defined for words $\w=x_1x_2\ldots x_n$ with
letters in $X$ by
\begin{align*}
{\inv}(\w)=\sum_{1\leq i<j\leq n}\chi(x_i>
x_j)\quad\text{and}\quad{\maj}(\w)&=\sum_{i=1}^{n-1}i.\chi(x_i>
x_{i+1}),
\end{align*}
where, as usual, $">"$ is the natural order on $X$ with
$r>r-1>\cdots>2>1$, and $\chi(A)=1$ if $A$ is true, and $\chi(A)=0$
otherwise.

 The major index, originally called \emph{greater index},
was introduced by MacMahon~\cite{Mac2}. As explained by Foata and
Krattenthaler (see \cite{FoKr} for a discussion), the origin of the
inversion number is not clear but probably MacMahon~\cite{Mac2,Ma}
was the first to consider inversions of words instead of just
permutations.

Let $\textbf{c}=(c(1),c(2),\ldots,c(r))$ be a sequence of $r$
non-negative integers and let $v$ be the non-decreasing word
$v=1^{c(1)}2^{c(2)}\ldots r^{c(r)}$. We will denote by $\R(v)$ (or
by $\R(\textbf{c})$ if there is no ambiguity) the
\emph{rearrangement class} of $v$, that is, the set of all words
that can be obtained by permuting the letters of~$v$. A well-known
result of MacMahon states that the major index and the inversion
number are equidistributed (i.e. have the same generating function)
on each rearrangement class $\R(\textbf{c})$. More precisely,
MacMahon showed that the generating function of the statistics
$\maj$ and $\inv$ on each $\R(\textbf{c})$ is given by
\begin{equation}\label{eq:Mac1} \sum_{\w\in\R(\c)}q^{\inv
(\w)}=\sum_{\w\in\R(\c)}q^{\maj (\w)}={c(1)+c(2)+\cdots+c(r)\brack
c(1),c(2),\cdots,c(r)}_{q}
\end{equation}
where, as usual in $q$-theory, the $q$-multinomial coefficient is
given by
$${n_1+n_2+\cdots+n_k\brack
n_1,n_2,\ldots,n_k}_{q}=\frac{[n_1+n_2+\ldots+n_k]_q!}{[n_1]_q![n_2]_q!\cdots
[n_k]_q!},$$ and the $q$-factorial $[n]_q!$  by
$[n]_q!:=(1+q)(1+q+q^2)\cdots (1+q+q^2+\cdots +q^{n-1})$. In honor
of MacMahon, a statistic which is equidistributed with $\inv$ (or
$\maj$) on each $\R(\textbf{c})$ is said to be \emph{mahonian}.

 In  1996,  Foata and Zeilberger~\cite{FoZe} introduced natural
generalizations of both "$\inv$" and "$\maj$", parametrized by
relations, as follows. Recall that a \emph{relation} $U$ on $X$ is a
subset of the cartesian product $X\times X$. For $a,b\in X$, if we
have $(a,b)\in U$, we say that $a$ is in relation $U$ to $b$, and we
express this also by $a\U b$. For each such relation $U$, then
associate the following statistics defined on each word
$\w=x_1\ldots x_n$ by
\begin{align*}
{\inv}'_{U}(\w)=\sum_{1\leq i<j\leq n}\chi(x_i\U x_j)
\quad\text{and}\quad {\maj}'_{U}(\w)&=\sum_{i=1}^{n-1}i.\chi(x_i\U
x_{i+1}).
\end{align*}
 For instance, where $U=">"$ is the natural order on $X$,
then ${\maj}'_{>}=\maj$ and ${\inv}'_{>}=\inv$. The statistics
${\maj}'_{U}$ and ${\inv}'_{U}$ are called \emph{graphical major
index} and \emph{graphical inversion number} since a relation on $X$
can be represented by a directed graph on $X$.

MacMahon's result \eqref{eq:Mac1} motivates Foata and Zeilberger
\cite{FoZe} to pose the following question:

\emph{For which relations $U$ on $X$ the statistics $\maj'_U$ and
$\inv'_U$ are equidistributed on each rearrangement class $\R(\c)$?}

Generalizing MacMahon's result, they have fully characterized such
relations. In order to present their result, we first recall the
following definition due to Foata and Zeilberger~\cite{FoZe}.

\begin{defn}\label{defn:bipartitional}
A relation $U$ on $X$ is said to be \emph{bipartitional} if there
exists an ordered partition $(B_1,B_2,\ldots,B_k)$ of $X$ into
blocks $B_l$ together with a sequence
$(\beta_1,\beta_2,\ldots,\beta_k)$ of 0's and 1's such that $x\U y$
if and only either (1) $x\in B_l$, $y\in B_{l'}$ and $l<l'$, or (2)
$x,y\in B_l$ and $\beta_l=1$.
\end{defn}

In this paper, we will use the following axiomatic characterization
of \emph{bipartitional relations} due to Han~\cite{Han}.

\begin{prop}\label{prop:bipartitional} \emph{A relation $U$ on $X$ is bipartitional if and only if (1) it
is transitive, i.e. $x\U y$ and $y\U z$ imply $x\U z$, and (2) for
each $x,y,z\in X$, $x\U y$ and $z\nU y$ imply $x\U z$}.
\end{prop}

Then Foata and Zeilberger \cite[Theorem~2]{FoZe} proved the following.\\[-0.4cm]

\textbf{Theorem A.} \emph{Let $U$ be a relation on $X$. The
statistics ${\maj}'_U$ and ${\inv}'_U$ are equidistributed on each
rearrangement class $\R(\c)$ if and only if
$U$ is bipartitional.}\\[-0.3cm]

 In this paper, we are interesting with statistics which are obtained by
summing a graphical major index and a graphical inversion number. In
order to motivate this work, we present here two such statistics.
The first one is the \emph{Rawlings major index}. In \cite{Raw1},
Rawlings have introduced statistics, denoted $k$-$\maj$ ($k\geq1$),
which interpolate the major index and the inversion number and
defined for words $\w=x_1\cdots x_n$ with letters in $X$ by
\begin{align*} k\text{-}{\maj(\w)}=&\sum_{i=1}^{n-1}i.\chi(x_i\geq
x_{i+1}+k)+\sum_{1\leq i<j\leq n}\chi(x_j+k>x_i>x_{j}).
\end{align*}
Note that $1\text{-}{\maj}=\maj$ while $r\text{-}{\maj}=\inv$.  Now,
if we set
\begin{align*}
U_k=\{(x,y)\in X^2\,/\,x\geq y+k\}\quad\text{and}\quad
V_k=\{(x,y)\in X^2\,/\,y+k>x>y\},
\end{align*}
we have $k\text{-}{\maj}={\maj}'_{U_k}+{\inv}'_{V_k}$. In
\cite{Raw2}, Rawlings proved that for each integer $k\geq 1$,
$k$-$\maj$ is a mahonian statistic. Since $U_k\cup V_k$ is the
natural order "$>$" on $X$, Rawlings 's result can be rewritten
${\maj}'_{U_k}+{\inv}'_{V_k}$ and ${\inv}'_{U_k\cup V_k}$ are
equidistributed on each rearrangement~class.

The second statistic is more recent and defined on words with
letters in a different alphabet. Let
$\mathcal{A}=\{A_1,A_2,\cdots,A_r\}$ be a collection of non-empty,
finite and mutually disjoints sets of non-negative integers.
Combining two statistics introduced by Steingrimsson \cite{Stein},
Zeng and the author~\cite{KZ} have defined a statistic, denoted
$\MAJ$, on words $\pi=B_1 B_2\cdots B_k$ with letters in
$\mathcal{A}$ by {\small$$\MAJ(\pi)=\sum_{1\leq i\leq
k-1}i.\chi(\min(B_i)>\max(B_{i+1}))+
 \sum_{1\leq i<j\leq k}\chi(\max(B_j)\geq\min(B_i)>\min(B_{j})).$$}
For instance, if $\pi=\{3,9\}\,\{2\}\,\{1,4,8\}\,\{7\}\,\{5,6\}$,
then $\MAJ(\pi)=(1+4)+(2)=7$. Let $U_{\mathcal{A}}$ and
$V_{\mathcal{A}}$ be the relations defined on $\mathcal{A}$ by
\begin{align*}
(B,B')\in U_{\mathcal{A}}&\Leftrightarrow \min(B)> \max(B'),\\
 (B,B')\in V_{\mathcal{A}}&\Leftrightarrow \max(B')\geq\min(B)>\min(B').
\end{align*}
 Then we have $\MAJ={\maj}'_{U_{\mathcal{A}}}+{\inv}'_{V_{\mathcal{A}}}$. It was proved in
\cite[Theorem 3.5]{KZ} that
\begin{equation}\label{eq:KZ1}
\sum_{\pi\in\,\R(A_1A_2\cdots A_r)}q^{\MAJ(\pi)}=[r]_q!.
\end{equation}
Since $U_{\mathcal{A}}\cup V_{\mathcal{A}}$ is a total order on
$\mathcal{A}$, it follows from \eqref{eq:Mac1} that the generating
function of $\inv'_{U_{\mathcal{A}}\cup V_{\mathcal{A}}}$ on
$\R(A_1A_2\cdots A_r)$ is also given by the right-hand side of the
above identity. It is then natural to ask if
${\maj}'_{U_{\mathcal{A}}}+{\inv}'_{V_{\mathcal{A}}}$ and
${\inv}'_{U_{\mathcal{A}}\cup V_{\mathcal{A}}}$ are equidistributed
on each rearrangement class $\R(\w)$ for words $\w$
with letters in~$\mathcal{A}$.\\

In view of the above two examples, it is natural to ask: \emph{For
which relations $U$ and $V$ on $X$ the statistics
${\maj}'_{U}+{\inv}'_V$ and ${\inv}'_{U\cup V}$ are
\begin{itemize}
\item equidistributed on each rearrangement class $\R(\c)$?
\item mahonian?
\end{itemize}}

The purpose of this paper is to answer these questions by fully
characterizing all such relations $U$ and $V$ on $X$.\\

\subsection{Main results}
Denote by $X^*$ the set of all words with letters in $X$. In order
to simplify the readability of the paper, we introduce the following
definition.

\begin{defn}\label{def:maj-inv stat}
A statistic $\stat$ on $X^*$ is a \emph{maj-inv statistic} if there
exist two relations $U$ and $V$ on $X$ such that
$\stat={\maj}'_U+{\inv}'_V.$
\end{defn}
Clearly, the statistics $\inv$, $\maj$ and $k$-$\maj$ are maj-inv
statistics on $X^*$, while $\MAJ$ is a maj-inv statistic on
$\mathcal{A}^*$. In this paper, a kind of relations on $X$ have a
great interest for us. We call them the $\kappa$-\emph{extensible
relations}.

\begin{defn}
A relation $U$ on $X$ is said to be \emph{$\kappa$-extensible} if
there exists a relation $S$ on $X$ such that (1) $U\subseteq S$ and
(2) for any $x,y,z\in X$, $x \U y$ and $z\nU y$ $\implies$ $x\S z$
and $z\nS x$.

If a relation $S$ on $X$ satisfies conditions (1) and (2), we say
that $S$ is a \emph{$\kappa$-extension} of $U$ on $X$.
\end{defn}

We give here some examples of $\kappa$-extensible relations.\\[-0.6cm]
\begin{exam}\label{exam:kappa extension}

(a) Suppose $X=\{x,y,z\}$ and $U=\{(x,y)\}$. Then,
$S=\{(x,y),(x,z)\}$ is a $\kappa$-extension of $U$ on $X$.

(b) The natural order "$>$" is a $\kappa$-extension of the relation
$U_k=\{(x,y)\in X^2\,/\,x\geq y+k\}$ on $X$ for any $k>0$.

(c) Let $\mathcal{A}=\{A_1,A_2,\cdots,A_r\}$ be a collection of
non-empty and finite subsets of non-negative integers, and let
$U_{\mathcal{A}}$ and $S_{\mathcal{A}}$ be the relations on
$\mathcal{A}$ defined by $(B,B')\in U_{\mathcal{A}} \Leftrightarrow
\min(B)>\max(B')$ and $(B,B')\in S_{\mathcal{A}} \Leftrightarrow
\min(B)>\min(B')$. Then one can check that $S_{\mathcal{A}}$ is a
$\kappa$-extension of $U_{\mathcal{A}}$ on $\mathcal{A}$.

(d) Every total order is a $\kappa$-extension of itself.

\end{exam}

In fact the notion of $\kappa$-extensible relation can be viewed, by
means of the following result, as a generalization of the notion of
bipartitional relation.

\begin{prop}\label{prop:bip-ext}
A relation $U$ on $X$ is bipartitional if and only if it is a
$\kappa$-extension of itself.
\end{prop}

\pf Using Proposition~\ref{prop:bipartitional}, it suffices to see
that a relation $U$ is transitive if and only if for any $x,y,z\in
X$, $x\U y$ and $z\nU y$ imply $z\nU x$.
 Suppose $U$ is transitive and let $x,y,z$ satisfying $x\U y$ and $z\nU y$. Suppose
$z\U x$, then since $x\U y$, we have by transitivity $z\U y$ which
contradict $z\nU y$. Thus $z\nU x$. Reversely, suppose that $x\U y$
and $z\nU y$ imply $z\nU x$ for each $x,y,z$. Let $x_1,x_2,x_3$
verifying $x_1\U x_2$ and $x_2\U x_3$. Suppose $x_1\nU x_3$. Since
$x_2\U x_3$, it then follows that $x_1\nU x_2$ which is impossible.
Thus $x_1\U x_3$ and $U$ is transitive. \qed\\

We can now present the key result of the paper, which is a
generalization of Theorem~A.\\[-0.6cm]
\begin{thm}\label{thm:majinv}
Let $U$ and $S$ be two relations on $X$. The following conditions
are equivalent.
\begin{itemize}
\item[(i)]
 The statistics $ {\maj}'_{\,U}+{\inv}'_{S\setminus U}$ and
${\inv}'_S$ are equidistributed on each rearrangement class
$\R(\c)$.
\item[(ii)]$S$ is a $\kappa$-extension of $U$.
\end{itemize}
\end{thm}

  Let $U$ and $V$ be two non-disjoint relations on $X$ and let $(x,y)\in U\cap
V$. By definition, $({\maj}'_{\,U}+{\inv}'_{V})(xy)=1+1=2>1\geq
\inv'_{U\cup V}(x_1x_2)$ for any $x_1,x_2\in X$. It follows that if
$U\cap V\neq \emptyset$, the statistics ${\maj}'_{\,U}+{\inv}'_{V}$
and $\inv'_{U\cup V}$ are not equidistributed on $\R(xy)$. We then
obtain immediately from Theorem~\ref{thm:majinv} the following
result.

\begin{thm}\label{thm:majinv-a}
Let $U$ and $V$ be two relations on $X$. The following conditions
are equivalent.
\begin{itemize}
\item[(i)]The statistics ${\maj}'_{\,U}+{\inv}'_{V}$ and ${\inv}'_{U\cup
V}$  are equidistributed on each rearrangement class $\R(\c)$.
\item[(ii)] $U\cap V=\emptyset$ and $\U\cup V$ is a
$\kappa$-extension of $U$.
\end{itemize}
\end{thm}

Next, by noting that for a relation $S$ on $X$, the graphical
inversion number  $\inv'_{S}$ is mahonian if and only if $S$ is a
total order on $X$, we have obtained the following characterization
of mahonian maj-inv statistics.

\begin{thm}[Classification of mahonian maj-inv statistics I]\label{thm:mahonian maj-inv}
 The mahonian maj-inv statistics on $X^*$ are
exactly those which can be written
 $\maj'_{U}+\inv'_{S\setminus U}$, where $U$ and $S$ satisfy the following conditions:
\begin{itemize}
\item $S$ is a total order on $X$,
\item $S$ is a $\kappa$-extension of $U$.
\end{itemize}

Moreover, two mahonian maj-inv statistics $\maj'_{U}+
\inv'_{S\setminus U}$ and $\maj'_{V}+\inv'_{T\setminus V}$ are equal
on $X^*$ if and only if $S=T$ and $U=V$.

\end{thm}

\begin{exam}\label{exam:mahonian stat}
(a) It follows from Example~\ref{exam:kappa extension}(b) and the
above theorem that the statistics $k$-$\maj$, $k\geq1$, are
mahonian, which was first proved by Rawlings~\cite{Raw2}.

 (b) Let $\mathcal{A}=\{A_1,A_2,\cdots,A_r\}$ be a collection of
nonempty and finite subsets of non-negative integers, and let
$U_{\mathcal{A}}$ and $S_{\mathcal{A}}$ be the relations on
$\mathcal{A}$ defined as in Example~\ref{exam:kappa extension}(c).
It then follows from the above theorem and Example~\ref{exam:kappa
extension}(c) that $\MAJ$ is mahonian on $\mathcal{A}^*$, which is a
generalization of~\eqref{eq:KZ1}.
\end{exam}

In fact, we have obtained more precise results on mahonian maj-inv
statistics on $X^*$. Indeed, given a total order $S$ on $X$, we have
characterized all $\kappa$-extensible relations $U$ such that $S$ is
a $\kappa$-extension of $U$(see Proposition~\ref{prop: preuve de
maj-inv countable}). As consequence, we have obtained the following
result.

\begin{thm}[Classification of mahonian maj-inv statistics II]\label{thm:maj-inv countable}
The mahonian maj-inv statistics on $X^*$ are exactly the statistics
$\stat_{f,\,g}$ defined for words $\w=x_1\cdots x_n\in X^*$ by
\begin{align*}
{\stat}_{f,\,g}(\w)=\sum_{i=1}^{n-1}i.\chi(\,f(x_i)\geq
g(f(x_{i+1}))\,)+\sum_{1\leq i<j\leq
n}\chi(\,g(f(x_j))>f(x_i)>f(x_j)\,),
\end{align*}
with $f$ a permutation of $X$ and  $g:X\mapsto X\cup\{\infty\}$ a
map satisfying $g(y)>y$ for each $y\in X$.
\end{thm}

Taking $f=Id$, where $Id$ is the identity permutation, we obtain the
following.

\begin{cor}\label{cor:mah maj-inv id}
The statistics $\stat_{g}$ defined for $\w=x_1\cdots x_n\in X^*$ by
\begin{align*}
{\stat}_{g}(\w)=\sum_{i=1}^{n-1}i.\chi(\,x_i\geq
g(x_{i+1})\,)+\sum_{1\leq i<j\leq n}\chi(\,g(x_j)>x_i>x_j\,),
\end{align*}
with $g:X\mapsto X\cup\{\infty\}$ satisfying $g(y)>y$ for each $y\in
X$, are mahonian.
\end{cor}

For instance, the Rawlings major index $k$-$\maj$ is obtained by
taking in the previous result $g:X\mapsto X\cup\{\infty\}$ defined
by $g(x)=x+k$ if $x+k\leq r$ and $g(x)=\infty$ otherwise.

It is then easy to enumerate the mahonian maj-inv statistics on
$X^*$. Since there are exactly $|X|!$ maps $g:X\mapsto
X\cup\{\infty\}$ satisfying $g(y)>y$, we have the following result.

\begin{cor}\label{cor:nbre de mahonian stat}
For each total order $S$ on $X$, there are exactly $|X|!$ mahonian
maj-inv statistics on $X^*$ which can be written
$\maj'_{U}+\inv'_{S\setminus U}$ .
\end{cor}

\vspace{0.3cm}

The paper is organized as follows. In section~2 and section~~3, we
prove Theorem~\ref{thm:majinv}. In section~4, we prove
Theorem~\ref{thm:mahonian maj-inv}. In section~5, we characterize
all $\kappa$-extensible relations on $X$ and  prove
Theorem~\ref{thm:maj-inv countable} in section~6. Finally, in
section~7, we apply the results of this paper to give new original
mahonian statistics on permutations and words.

\begin{rmk}
As pointed by an anonymous referee, some proofs ( for instance the
proof of the "only if" part of Theorem~\ref{thm:majinv}) presented
in the paper have "simpler proofs" by using a computer algebra
system (see e.g. \cite{Han1}).
\end{rmk}

%%%%****%%%%****%%%%****%%%%****%%%%****%%%%****%%%%****%%%%*
%%%%
%%%% New Section
%%%%
%%%%****%%%%****%%%%****%%%%****%%%%****%%%%****%%%%****%%%%*

\section{Proof of the 'if' part of Theorem~\ref{thm:majinv}}

 The first direct combinatorial proof of MacMahon's
result on the equidistribution of the statistics $\maj$ and $\inv$,
that is a bijection which sends each word to another one in such a
way that the major index of the image equals the number of
inversions of the original, is due to Foata \cite{Fo}.

Let $U$ be a $\kappa$-extensible relation on $X$. In this section,
we adapt Foata's map, also called \emph{second fundamental
transformation} (see e.g. \cite{Lo}), to construct a bijection
$\Psi^U$ of each rearrangement class onto itself such that for each
$\kappa$-extension $S$ of $U$, we have
\begin{equation}
{\inv}'_S(\Psi^U(\w))=({\maj}'_U+{\inv}'_{S\setminus U})(\w).
\end{equation}

\subsection{Notations}

  The \emph{length} of a word $\w\in X^*$, denoted by $\lambda(\w)$, is its number of
letters. By convention, there is an unique word of length 0, the
\emph{empty word} $\epsilon$. If $Y$ and $Z$ are subsets of $X^*$,
we designate by $YZ$ the set of words $\w=\w'\w''$ with $\w'\in Y$
and $\w''\in Z$.

Each $x\in X$ determines a partition of $X$ in two subsets $L_x$ and
$R_x$ as follows: the set $R_x$ is formed with all $y\in X$ such
that $y\U x$ , while the set $L_x$ is formed with all $y\in X$ such
that $y\nU x$.

\subsection{The map $\Psi^U$}
 Let $\w$ be a word in
$X^*$ and $x\in X$. If $\w=\epsilon$, we set
$\gamma_{x}^U(\w)=\epsilon$. Otherwise two cases are to be
considered:
\begin{itemize}
\item[(i)] the last letter of $\w$ is in $R_x$,
\item[(ii)] the last letter of $\w$ is in $L_x$.
\end{itemize}

Let $(\w_1x_1,\w_2x_2,\ldots,\w_hx_h)$ be the factorization of $\w$
having the following properties:

 \begin{itemize}
\item In case (i) $x_1,x_2,\ldots,x_h$ are in $R_x$ and
$\w_1,\w_2,\ldots,\w_h$ are words in $L_x^*$.
\item In case (ii) $x_1,x_2,\ldots,x_h$ are in $L_x$ and
$\w_1,\w_2,\ldots,\w_h$ are words in $R_x^*$.
\end{itemize}
Call $x$-\emph{factorization} the above factorization. Clearly, each
word has an unique $x$-factorization. In both cases we have
$\w=\w_1x_1\w_2x_2\ldots \w_sx_s$, then define
$$
\gamma_x^U(\w)=x_1\w_1x_2\w_2\ldots x_s\w_s.
$$

 The map $\Psi^U$ is then defined by induction on the length of words
in the following way:
\begin{align}
\Psi^U(\epsilon)&=\epsilon \;, \label{eq:defbij1}\\
\Psi^U(\w x)&=\gamma_x^U(\Psi^U(\w))\,x\quad\text{for all $x\in X$
and $\w\in X^*$}.\label{eq:defbij2}
\end{align}

Note that Foata's map correspond to the case $U:=">"$ is the natural
order.
\begin{thm}\label{thm:Psi}
The map $\Psi^U$ is a bijection of $X^*$ onto itself such that for
each $\w\in X^*$, we have $\Psi^U(\w)\in\R(\w)$, both $\w$ and
$\Psi^U(\w)$ end with the same letter and for each
$\kappa$-extension $S$ of $U$, we have
\begin{equation}\label{eq:Psi}
{\inv}'_S(\Psi^U(\w))=({\maj}'_U+{\inv}'_{S\setminus U})(\w).
\end{equation}
\end{thm}

The proof of the above theorem is very similar to the proof in
\cite{Fo,Lo}. It is based on the following lemma. Let $S$ be a
$\kappa$-extension of $U$.  For each $\w=x_1\ldots x_n\in X^*$,
denote by $l_x(\w)$ (resp. $r_x(\w)$) the number of subscripts $j$
for which $x_j\in L_x$ (resp. $x_j\in R_x$) and $t_x(\w)$ designate
the number of subscripts $j$ such that $x_j\nU x$ and $x_j\S x$.
Note that we always have $l_x(\w)+r_x(\w)=\lambda(\w)$ and
$r_x(\w)+t_x(\w)$ is the number of subscripts $j$ for which $x_j\S
x$.
\begin{lem}\label{lem:gamma prop}
For each $\w\in X^*$ and $x\in X$, the following identities hold:
\begin{eqnarray}
  {\inv}'_S(\w x) &=& {\inv}'_S(\w)+r_x(\w)+t_x(\w), \label{eq:prop1}\\
  {\inv}'_S(\gamma^U_x(\w)) &=& {\inv}'_S(\w)-r_x(\w)\qquad\text{if}\;\w\in
  X^*L_x, \label{eq:prop2}\\
  {\inv}'_S(\gamma^U_x(\w)) &=& {\inv}'_S(\w)+l_x(\w)\qquad\text{if}\;\w\in X^*R_x, \label{eq:prop3}\\
  ({\maj}'_U+{\inv}'_{S\setminus U})(\w x) &=& ({\maj}'_U+{\inv}'_{S\setminus U})(\w)+t_x(\w)\qquad\text{if}\;\w\in
  X^*L_x, \label{eq:prop4}\\
  ({\maj}'_U+{\inv}'_{S\setminus U})(\w x)&=& ({\maj}'_U+{\inv}'_{S\setminus U})(\w)+ t_x(\w)+\lambda(\w)\quad\text{if}\;\w\in X^*R_x.\label{eq:prop5}
\end{eqnarray}
\end{lem}

\pf  By definition, we have the following identities:
   \begin{eqnarray*}
  {\inv}'_{U}(\w x)&=&{\inv}'_{U}(\w)+r_x(\w)\\
{\inv}'_{S\setminus U}(\w x)&=&{\inv}'_{S\setminus U}(\w)+t_x(\w)\\
    {\maj}'_U(\w x) &=& {\maj}'_U (\w)\qquad\text{if}\;\w\in
  X^*L_x, \\
  {\maj}'_U(\w x)&=& {\maj}'_U(\w)+\lambda(\w)\quad\text{if}\;\w\in
  X^*R_x,
\end{eqnarray*}
from which we derive immediately \eqref{eq:prop4} and
\eqref{eq:prop5}. To obtain \eqref{eq:prop1}, it suffices to note
that $\inv'_S={\inv}'_U+{\inv}'_{S\setminus U}$ (since $U\cap
(S\setminus U)=\emptyset$ and $U\subseteq S$). It remains to prove
\eqref{eq:prop2} and \eqref{eq:prop3}.

Suppose $\w\in X^*L_x$ and let $(\w_1x_1,\w_2x_2,\ldots,\w_sx_s)$ be
the $x$-factorization of $\w$. First, assume that
 \begin{equation}\label{eq:lemmepreuve1}
 {\inv}'_S(x_i\w_i)={\inv}'_S(\w_ix_i)-\lambda(\w_i)\quad\text{for}\;1\leq i\leq s.
 \end{equation}
  Since $\gamma_x(\w)=x_1\w_1x_2\w_2\cdots x_h\w_h$,
 it is not hard to see that ${\inv}'_S(\gamma_x(\w))$ is equal to ${\inv}'_S(\w x)$
 decreased by $\lambda(\w_1)+\lambda(\w_2)+\cdots+\lambda(\w_s)$.
 Since $s=l_x(\w)$, we get
 \begin{align*}
{\inv}'_S(\gamma_x(\w))={\inv}'_S(\w x)-(\lambda(\w)-s)={\inv}'_S(\w
x)-r_x(\w),
 \end{align*}
which is exactly \eqref{eq:prop2}. We now prove
\eqref{eq:lemmepreuve1}. Let $\tau=\tau_1\tau_2\cdots\tau_m\in
R_x^*$ and $y\in L_x$.  By definition, we have $\tau_i \U x$ for
each $i$ and $y\nU x$. Since $S$ is a $\kappa$-extension of $U$, it
follows that for each $i$, $\tau_i\S y$ and $y\nS \tau_i$. We then
have ${\inv}'_S(y\tau)={\inv}'_S(\tau)$ and ${\inv}'_S(\tau
y)={\inv}'_S(\tau)+ m ={\inv}'_S(y\tau)+\lambda(\tau)$. Equation
\eqref{eq:lemmepreuve1} is obtained by noting that in the
$x$-factorization of $\w\in X^*L_x$, the words $\w_1,\ldots,\w_h$
are in $R_x^*$ and the letters $x_1,\ldots,x_h$ are in $L_x$.

  Equation \eqref{eq:prop3} has an analogous proof.
Suppose $\w\in X^*R_x$ and let $(\w_1x_1,\w_2x_2,\ldots,\w_hx_h)$ be
the $x$-factorization of $\w$. First, assume that
 \begin{equation}\label{eq:lemmepreuve2}
 {\inv}'_S(x_i\w_i)={\inv}'_S(\w_ix_i)+\lambda(\w_i)\quad\text{for}\;1\leq i\leq h.
 \end{equation}
  Since $\gamma_x(\w)=x_1\w_1x_2\w_2\cdots x_h\w_h$,
 it is not hard to see that ${\inv}'_S(\gamma_x(\w))$ is equal to ${\inv}'_S(\w x)$
 increased by $\lambda(\w_1)+\lambda(\w_2)+\cdots+\lambda(\w_s)$.
 Since $h=r_x(\w)$, we get
 \begin{align*}
{\inv}'_S(\gamma^U_x(\w))={\inv}'_S(\w
x)-(\lambda(\w)-h)={\inv}'_S(\w x)+l_x(\w),
 \end{align*}
which is exactly \eqref{eq:prop3}. It then remains to prove
\eqref{eq:lemmepreuve2}. Let $\tau=\tau_1\tau_2\cdots\tau_m\in
L_x^*$ and $y\in R_x$. By definition, we have $y\U x$ and $\tau_i\nU
x$ for each $i$. Since $S$ is a $\kappa$-extension of $U$, it
follows that for each $i$, $y\S\tau_i$. It is then easy to obtain
${\inv}'_S(\tau y)={\inv}'_S(\tau)$ and
${\inv}'_S(y\tau)={\inv}'_S(\tau)+m={\inv}'_S(\tau y)+m$. Equation
\eqref{eq:lemmepreuve2} is obtained by noting that in the
$x$-factorization of $\w\in X^*R_x$, the words $\w_1,\ldots,\w_h$
are in $L_x^*$ and the letters $x_1,\ldots,x_h$ are in $R_x$.

\qed

\vspace{0.3cm}

\emph{Proof of Theorem~\ref{thm:Psi}:}  By construction, both $\w$
and $\Psi^U(\w)$ end with the same letter. Let $X_n$ be the set of
words in $X^*$ with length $n$. It is sufficient to verify by
induction on $n$ that for all $n\geq0$, the restriction $\Psi^U_n$
of $\Psi^U$ to $X_n$ is a permutation of $X_n$
 satisfying: for any $\w\in X_n$,
\begin{align}
\Psi^U_n(\w)\in \R(\w)\quad\text{and}\quad
{\inv}'_S(\Psi^U_n(\w))=({\maj}'_U+{\inv}'_{S\setminus
U})(\w).\label{eq:preuvepsi1}
\end{align}

Since the induction is based on Lemma~\ref{lem:gamma prop} and is
very similar to the proof concerning the second fundamental
transformation, we refer the reader to \cite{Fo,Lo}.
 \qed

%%%%****%%%%****%%%%****%%%%****%%%%****%%%%****%%%%****%%%%*
%%%%
%%%% New Section
%%%%
%%%%****%%%%****%%%%****%%%%****%%%%****%%%%****%%%%****%%%%*

\section{Proof of the 'only if' part of Theorem~\ref{thm:majinv}}

 Let $U$ and $S$ be two relations on $X$ such that the statistics
${\maj}'_U+{\inv}'_{S\setminus U}$ and ${\inv}_S$ are
equidistributed on each rearrangement class $\R(\w)$, $\w\in X^*$.
We prove here that this imply  that $S$ is a $\kappa$-extension of
$U$.

\subsection{The relation $U$ is contained in $S$}
By definition of the graphical statistics, we have that for all word
$\w$ of length~2, ${\maj}'_U(\w)={\inv}'_{U}(\w)$. Moreover, for
each pair $(A,B)$ of disjoints relations, we have
${\inv}'_A+{\inv}'_{B}={\inv}'_{A\cup B}$. It then follows that for
all $\w\in X^*$, $\lambda(\w)=2$,
\begin{align*}
({\maj}'_U+{\inv}'_{S\setminus
U})(\w)=({\inv}'_U+{\inv}'_{S\setminus U})(\w)={\inv}'_{S\cup
U}(\w)=({\inv}'_{S}+{\inv}'_{U\setminus S})(\w).
\end{align*}
The equidistribution of ${\maj}'_U+{\inv}'_{S\setminus U}$ and
$\inv'_S$ on each $\R(\w)$, $\lambda(\w)=2$, then implies that
${\inv}'_{U\setminus S}(\w)=0$ for all $\w\in X^*$, $\lambda(\w)=2$,
and thus, $U\setminus S=\emptyset$, i.e., $U\subseteq S$.

\subsection{For any $x,y,z\in X$, $x\U y$ and $z\nU y$ imply $x\S z$ and $z\nS x$}

 To simplify the readability of the rest of the proof, we set
$V:=S\setminus U$, i.e. $U\cap V=\emptyset$ and $U\cup V=S$. In
particular, for any $x_1,x_2\in X$,
\begin{align}\chi(x_1\S
x_2)=\chi(x_1\U x_2) +\chi(x_1\V x_2)\quad\text{and}\quad\chi(x_1\U
x_2).\chi(x_1\V x_2)=0.\label{eq:simplifications}\end{align}

 Let $x,y,z\in X$ verifying $x\U y$ and $z\nU y$. First,
note that $x$ and $z$ are distinct, otherwise we have $x\U y$ and
$x\nU y$. Thus $x\neq z$.

\subsubsection{The case $x=y$} We then have $x\U x$ and $z\nU x$ and thus\\
$({\maj}'_U+{\inv}'_{S\setminus U})(zxx)=\chi(z\U x)+2\chi(x\U x)+2
\chi(z\V x)+\chi(x\V x)=2+2 \chi(z\V x)$. Since $\inv'_S(\w)\leq 3$
for each word $\w$ of length $3$, it follows that $\chi(z\V x)=0$,
i.e. $z\nV x$. But $z\nU x$ and thus $z\nS x$. Now, suppose $x\nS
z$. It follows that
$\inv'_S(xxz)=\inv'_S(xzx)=\inv'_S(zxx)=1<2=({\maj}'_U+{\inv}'_{S\setminus
U})(zxx)$, which contradict the equidistribution of our two
statistics on $\R(x^2z)$. Thus we have $x\S z$ and $z\nS x$ as
desired.

\subsubsection{The case $x\neq y$} Two cases are
to be considered.

\emph{Suppose $y=z$}. We then have $x\U z$ and $z\nU z$. Since
$U\subseteq S$, we have $x\S z$. It then suffices to show that $z\nS
x$. Suppose $z\S x$. We then have
\begin{align*}({\maj}'_U+{\inv}'_{S\setminus U})(zxz)&=\chi(z\U x)+2\chi(x\U
z)+\chi(z\V x)+\chi(z\V z)+\chi(x\V z)\\&=2+\chi(z\U x)+\chi(z\V
x)+\chi(z\V z) =2+\chi(z\S x)+\chi(z\V z)\\&=3+\chi(z\V z),
\end{align*}
and thus $z\nS z$. Then, it is not hard to see that this imply that
$\inv'_S\leq 2$ on $\R(xz^2)$ which, considering
$({\maj}'_U+{\inv}'_{S\setminus U})(zxz)=~3$, contradict the
equidistribution of our statistics on $\R(xz^2)$. It follows that
$z\nS x$ as desired. It then remains to consider
the last case.\\[-0.3cm]

 \emph{Suppose $y\neq z$}. Then \emph{$x,y,z$ are three distinct elements} satisfying $x\U y$ and $z\nU y$.
The next table gives the distribution of
${\maj}'_U+{\inv}'_{S\setminus U}$ and $\inv'_S$
on $\R(xyz)$ after some simplifications obtained by using \eqref{eq:simplifications}.\\

\centerline{%
\vbox{\halign{\vrule\strut\ \hfil$#$\hfil\ \vrule &\kern 11pt
\hfil$#$\hfil\kern 10pt \vrule &\kern 11pt \hfil$#$\hfil\kern 10pt
\vrule  \cr \noalign{\hrule}
 \w&({\maj}'_U+{\inv}'_{S\setminus
U})(\w)&\inv'_S(\w)\cr \noalign{\hrule}
 xyz&1+\chi(y\S z)+\chi(y \U z)+\chi(x \V z) & 1+\chi(x\S z)+\chi(y \S z) \cr \noalign{\hrule}
 xzy&\chi(x\S z)+\chi(z \V y) & 1+\chi(x\S z)+\chi(z \V y)\cr\noalign{\hrule}
 yxz&\chi(y\S x)+\chi(x \S z)+\chi(x \U z)+\chi(y \V z)& \chi(y\S x)+\chi(y \S z)+\chi(x \S z)\cr \noalign{\hrule}
 yzx&\chi(y\S z)+\chi(z \S x)+\chi(z \U x)+\chi(y \V x)& \chi(y\S z)+\chi(y \S x)+\chi(z \S x)\cr \noalign{\hrule}
 zxy&2+\chi(z\S x)+\chi(z \V y)&1+\chi(z\S x)+\chi(z \V y)\cr\noalign{\hrule}
 zyx&\chi(y\S x)+\chi(y \U x)+\chi(z \V y)+\chi(z \V x)&\chi(z\V y)+\chi(z \S x)+\chi(y\S x)\cr \noalign{\hrule} }}}

\vspace{0.5cm}

 (a) Suppose $x\S z$ and $z\S x$. We then have
 $({\maj}'_U+{\inv}'_{S\setminus
U})(zxy)=3+\chi(z \V y)$ and since $\inv'_S\leq 3$ on $\R(xyz)$, we
have $z\nV y$ and thus $z\nS y$. Using identities $x\S z$, $z\S x$
and $z\nS
y$, we obtain the following table\\

\centerline{%
\vbox{\halign{\vrule\strut\ \hfil$#$\hfil\ \vrule &\kern 11pt
\hfil$#$\hfil\kern 10pt \vrule &\kern 11pt \hfil$#$\hfil\kern 10pt
\vrule  \cr \noalign{\hrule}
 \w&({\maj}'_U+{\inv}'_{S\setminus
U})(\w)&\inv'_S(\w)\cr \noalign{\hrule}
 xyz&1+\chi(y\S z)+\chi(y \U z)+\chi(x \V z) & 2+\chi(y \S z) \cr \noalign{\hrule}
 xzy&1& 2\cr\noalign{\hrule}
 yxz&1+\chi(y\S x)+\chi(x \U z)+\chi(y \V z)& 1+\chi(y\S x)+\chi(y \S z)\cr \noalign{\hrule}
 yzx&1+\chi(y\S z)+\chi(z \U x)+\chi(y \V x)& 1+\chi(y \S x)+\chi(y\S z)\cr \noalign{\hrule}
 zxy&3&2\cr\noalign{\hrule}
 zyx&\chi(y\S x)+\chi(y \U x)+\chi(z \V x)&1+\chi(y\S x)\cr \noalign{\hrule} }}}

\vspace{0.2cm}

which imply that $y\nS x$ (otherwise, $\inv'_S\geq 2$ on
 $\R(xyz)$ and $({\maj}'_U+{\inv}'_{S\setminus
U})(xzy)=1$, which is impossible) and thus, by using $y\nS x$, we get\\

\centerline{%
\vbox{\halign{\vrule\strut\ \hfil$#$\hfil\ \vrule &\kern 11pt
\hfil$#$\hfil\kern 10pt \vrule &\kern 11pt \hfil$#$\hfil\kern 10pt
\vrule  \cr \noalign{\hrule}
 \w&({\maj}'_U+{\inv}'_{S\setminus
U})(\w)&\inv'_S(\w)\cr\noalign{\hrule}
 xyz&1+\chi(y\S z)+\chi(y \U z)+\chi(x \V z) & 2+\chi(y \S z) \cr \noalign{\hrule}
 xzy&1& 2\cr\noalign{\hrule}
 yxz&1+\chi(x \U z)+\chi(y \V z)& 1+\chi(y \S z)\cr \noalign{\hrule}
 yzx&1+\chi(y\S z)+\chi(z \U x)+\chi(y \V x)& 1+\chi(y\S z)\cr \noalign{\hrule}
 zxy&3&2\cr\noalign{\hrule}
 zyx&\chi(z \V x)&1\cr \noalign{\hrule} }}}

\vspace{0.2cm}

 Since $({\maj}'_U+{\inv}'_{S\setminus U})(zxy)=3$, we have by equidistribution
 of our two statistics, $y\S z$. It follows that
 $zyx$ is the unique world in $\R(xyz)$ for which $\inv'_S(zyx)=1$,
 while $({\maj}'_U+{\inv}'_{S\setminus U})(zyx)\leq({\maj}'_U+{\inv}'_{S\setminus U})(xzy)\leq 1$,
 which contradict the equidistribution of our two statistics on $\R(xyz)$.\\

(b) Suppose $x\nS z$ and $z\S x$. By a similar reasoning than in
(a), we have $z\nV y$ and
thus $z\nS y$, which lead to the following table.\\

\centerline{%
\vbox{\halign{\vrule\strut\ \hfil$#$\hfil\ \vrule &\kern 11pt
\hfil$#$\hfil\kern 10pt \vrule &\kern 11pt \hfil$#$\hfil\kern 10pt
\vrule  \cr \noalign{\hrule}
 \w&({\maj}'_U+{\inv}'_{S\setminus
U})(\w)&\inv'_S(\w)\cr \noalign{\hrule}
 xyz&1+\chi(y\S z)+\chi(y \U z)& 1+\chi(y \S z) \cr \noalign{\hrule}
 xzy&0& 1\cr\noalign{\hrule}
 yxz&\chi(y\S x)+\chi(y \V z)& \chi(y\S x)+\chi(y \S z)\cr \noalign{\hrule}
 yzx&1+\chi(y\S z)+\chi(z \U x)+\chi(y \V x)& 1+\chi(y \S x)+\chi(y\S z)\cr \noalign{\hrule}
 zxy&3&2\cr\noalign{\hrule}
 zyx&\chi(y\S x)+\chi(y \U x)+\chi(z \V x)&1+\chi(y\S x)\cr \noalign{\hrule} }}}

\vspace{0.2cm}

 Since $({\maj}'_U+{\inv}'_{S\setminus
U})(zxy)=3$, we must have $y\S z$ and $y\S x$, which imply that
$\inv'_S\geq 1$ on $\R(xyz)$, which is impossible since
$(\maj'_U+\inv'_{S\setminus
U})(xzy)=0$.\\

 (c) Suppose $x\nS z$ and $z\nS x$. We then get the following table.\\

  \centerline{%
\vbox{\halign{\vrule\strut\ \hfil$#$\hfil\ \vrule &\kern 11pt
\hfil$#$\hfil\kern 10pt \vrule &\kern 11pt \hfil$#$\hfil\kern 10pt
\vrule  \cr \noalign{\hrule}
 \w&({\maj}'_U+{\inv}'_{S\setminus
U})(\w)&\inv'_S(\w)\cr\noalign{\hrule}
 xyz&1+\chi(y\S z)+\chi(y \U z)& 1+\chi(y \S z) \cr \noalign{\hrule}
 xzy&\chi(z \V y) & 1+\chi(z \V y)\cr\noalign{\hrule}
 yxz&\chi(y\S x)+\chi(y \V z)& \chi(y\S x)+\chi(y \S z)\cr \noalign{\hrule}
 yzx&\chi(y\S z)+\chi(y \V x)& \chi(y \S x)+\chi(y\S z)\cr \noalign{\hrule}
 zxy&2+\chi(z \V y)&1+\chi(z \V y)\cr\noalign{\hrule}
 zyx&\chi(y\S x)+\chi(y \U x)+\chi(z \V y)&\chi(z\V y)+\chi(y\S x)\cr \noalign{\hrule}
}}}

\vspace{0.2cm}

It then follows that $\inv'_S\leq 2$ on $\R(xyz)$, and thus, by
considering $({\maj}'_U+{\inv}'_{S\setminus U})(zxy)$ and
$({\maj}'_U+{\inv}'_{S\setminus U})(xyz)=1+2\chi(y\U z)+\chi(y\V
z)$, we must have $y\nU z$ and
$z\nV y$, which lead to the following table.\\

 \centerline{%
\vbox{\halign{\vrule\strut\ \hfil$#$\hfil\ \vrule &\kern 11pt
\hfil$#$\hfil\kern 10pt \vrule &\kern 11pt \hfil$#$\hfil\kern 10pt
\vrule  \cr \noalign{\hrule}
  \w&({\maj}'_U+{\inv}'_{S\setminus
U})(\w)&\inv'_S(\w)\cr\noalign{\hrule}
 xyz&1+\chi(y\V z)& 1+\chi(y \V z) \cr \noalign{\hrule}
 xzy&0 & 1\cr\noalign{\hrule}
 yxz&\chi(y\S x)+\chi(y \V z)& \chi(y\S x)+\chi(y \V z)\cr \noalign{\hrule}
 yzx&\chi(y\V z)+\chi(y \V x)& \chi(y \S x)+\chi(y\V z)\cr \noalign{\hrule}
 zxy&2&1\cr\noalign{\hrule}
 zyx&\chi(y\S x)+\chi(y \U x)&\chi(y\S x)\cr \noalign{\hrule} }}}

\vspace{0.2cm}

 Since $({\maj}'_U+{\inv}'_{S\setminus
U})(xzy)=0$, it follows that $y\nS x$, and thus $\inv'_S\leq 1$ on
$\R(xyz)$, which is impossible
 since $({\maj}'_U+{\inv}'_{S\setminus
U})(zxy)=2$.\\

 (d) Finally, we have $x\S z$ and $z\nS x$, and thus $S$ is a $\kappa$-extension of
$U$. This conclude the proof of the 'only if' part of
Theorem~\ref{thm:majinv}.

%%%%****%%%%****%%%%****%%%%****%%%%****%%%%****%%%%****%%%%*
%%%%
%%%% New Section
%%%%
%%%%****%%%%****%%%%****%%%%****%%%%****%%%%****%%%%****%%%%*

\section{mahonian maj-inv statistics}

 This section is dedicated to the proof of Theorem~\ref{thm:mahonian maj-inv}.
We begin with two lemmas.

\begin{lem}\label{lem:mahonian inv}
 Let $S$ be a relation on $X$. Then $\inv'_S$ is mahonian on $X^*$ if and only if
$S$ is a total order.
\end{lem}

\begin{lem}\label{lem:mahonian maj-inv}
 Let $U$ and $V$ be two relations on $X$. Suppose that the statistic $\maj'_U+\inv'_V$
is mahonian on $X^*$. Then, $U\cap V=\emptyset$, $S:=U\cup V$ is a
total order and a $\kappa$-extension of~$U$.
\end{lem}

It is now easy to prove the first part of Theorem~\ref{thm:mahonian
maj-inv}. Indeed, suppose that $S$ is a total order on $X$ and a
$\kappa$-extension of $U$. Then, it follows from
Theorem~\ref{thm:majinv} that ${\maj}'_U+\inv'_{S\setminus U}$ is
equidistributed with $\inv'_S$ which is mahonian by
Lemma~\ref{lem:mahonian inv}, and thus ${\maj}'_U+\inv'_{S\setminus
U}$ is mahonian as desired.  Reversely, suppose $\maj'_U+\inv'_V$ is
mahonian on $X^*$. We  then have by Lemma~\ref{lem:mahonian maj-inv}
$V=S\setminus U$ where $S:=U\cup V$ is a total order on $X$ and a
$\kappa$-extension of $U$.

We thus have proved that the mahonian maj-inv statistics on $X^*$
are exactly those which can be written
 $\maj'_{U}+\inv'_{S\setminus U}$, with $S$ a total order on $X$
and a $\kappa$-extension of $U$.\\

We now prove the second part of Theorem~\ref{thm:mahonian maj-inv}.
Let $S$ and $T$ be two total orders on $X$ and suppose $S$ (resp.
$T$) is a $\kappa$-extension of $U$ (resp. $V$). It suffices to show
that if $\maj'_{U}+\inv'_{S\setminus U}$ and
$\maj'_{V}+\inv'_{T\setminus V}$ are equal on $X^*$ then $S=T$ and
$U=V$.

 First suppose that $S\neq T$. Then we can assume without loss
of generality that there exist $x,y\in X$ such that $x\S y$ and
$x\nT y$. Since $U\subseteq S$ and $V\subseteq T$, we then have\\ $
\maj'_{U}+\inv'_{S\setminus U}(xy)=1\neq
0=\maj'_{V}+\inv'_{T\setminus V}(xy)$ and thus the two statistics
are different which contradict the hypothesis, thus $S=T$.

Suppose now $U\neq V$. Then we can assume without loss of generality
that there exist $x,y\in X$ such that $x\U y$ and $x\nV y$. Since
$S=T$ is a total order and an extension of $U$ and $V$ we also have
$x\S y$, $(x,y)\in S\setminus V$ and $y\nS y$. It follows
that\\
$ \maj'_{U}+\inv'_{S\setminus U}(xy^2)=1\neq
2=\maj'_{V}+\inv'_{T\setminus V}(xy^2), $ which is impossible thus
$U=V$, as desired.

In order to complete the proof of Theorem~\ref{thm:mahonian
maj-inv},
it then remains to prove the two above lemmas.\\

\emph{Proof of Lemma~\ref{lem:mahonian inv}.} It suffices to see
that $\inv'_S$ is mahonian imply that $S$ is a total order since the
reciprocal is an easy consequence of \eqref{eq:Mac1}. Suppose that
$\inv'_S$ is mahonian, i.e. for each $\textbf{c}$,
\begin{equation}\label{eq:invS-mah} \sum_{\w\in\R(\c)}q^{{\inv}'_S
(\w)}={c(1)+c(2)+\cdots+c(r)\brack c(1),c(2),\cdots,c(r)}_{q}.
\end{equation}

 Suppose there exist $x,y\in X$, $x\neq y$, such that $x\nS y$ and $y\nS x$. We then have
$\inv'_S(xy)=\inv'_S(yx)=0$, which contradict \eqref{eq:invS-mah}
(take $\w=xy$). Thus for each $x,y\in X$, we have $x\S y$ or $y\S
x$, i.e., $S$ is total.

 Suppose there exist $x\in
X$ such that $x\S x$, then $\inv'_S(xx)=1$, which contradict
\eqref{eq:invS-mah} (take $\w=x^2$). Thus $x\nS x$ and $S$ is
irreflexive.

 Suppose there exist $x,y\in
X$, $x\neq y$, such that $x\S y$ and $y \S x$. We then have
$\inv'_S(xy)=\inv'_S(yx)=1$, which contradict \eqref{eq:invS-mah}
(take $\w=xy$). Thus if $x\S y$ we have $y \nS x$, i.e. S is
antisymmetric.

  Let $x,y,z\in X$ satisfying $x\S y$ and $y\S z$.
Suppose $x\nS z$.  Since $S$ is irreflexive, we have $x\neq y$ and
$y\neq z$. Since $S$ is antisymmetric, we have $x\neq z$ (otherwise
we have $x\S y$ and $y\S x$). Then $x,y,z$ are distinct. We also
have $y\nS x$ and $z\nS y$ ($S$ is antisymmetric) and $z\S x$ ($S$
is total). After simple computations (we left the details to the
reader), we then get
\begin{align*}
\sum_{\w\in\R(xyz)}q^{\inv'_S(\w)}&=3q+3q^2\neq{3\brack
1,1,1}_{q}=1+2q+2q^2+q^3,
\end{align*}
which contradict \eqref{eq:invS-mah} (take $\w=xyz$). Thus $x\S z$
and $S$ is transitive.

\qed

\emph{Proof of Lemma~\ref{lem:mahonian maj-inv}.} Suppose $U\cap
V\neq \emptyset$ and let $(x,y)\in U\cap V$. We then have
$\maj'_U(xy)+\inv'_V(xy)=1+1=2$, which contradict
\eqref{eq:invS-mah} (take $\w=xy$ if $x\neq y$ and $\w=xx$ if $x=y$)
and thus, $U$ and $V$ are disjoint.

The proof of "$S$ is a total order on $X$" is essentially the same
than the proof of Lemma~\ref{lem:mahonian inv}, so we left the
details to the reader.

It then remains to show that $S=U\cup V$ is a $\kappa$-extension of
$U$. Since $S$ is total, it follows from Lemma~\ref{lem:mahonian
inv} that $\inv'_S$ are mahonian on $X^*$. Then by applying the part
"(i) imply (ii)" of Theorem~\ref{thm:majinv} to $\maj'_U+\inv'_V$
and $\inv'_S$, we obtain that $S$ is a $\kappa$-extension of~$U$.

\qed

%%%%****%%%%****%%%%****%%%%****%%%%****%%%%****%%%%****%%%%*
%%%%
%%%% New Section
%%%%
%%%%****%%%%****%%%%****%%%%****%%%%****%%%%****%%%%****%%%%*

\section{ $\kappa$-extensible relations}

Theorem~\ref{thm:majinv} and Theorem~\ref{thm:mahonian maj-inv}
motivate to pose the following question: \emph{When does a relation
have a $\kappa$-extension?}

 Suppose $X=\{x,y,z\}$ and consider the relation
$U=\{(x,y),(y,z)\}$ on $X$. Then, one can check by considering all
the relations on $X$ containing $U$ (there are $2^{3^2-2}=128$ such
relations) that $U$ has no $\kappa$-extension. In this part, we give
an axiomatic characterization of $\kappa$-extensible relations.

\begin{defn}
The $\kappa$-closure of a relation $U$ on a set $X$ is the relation
denoted by $cl_{\kappa}(U)$ and defined by
\begin{equation}
cl_{\kappa}(U):=U \cup \{(x,y) /\quad \text{$\exists\,z\in X$ such
that $x\U z$ and $y \nU z$}\}.
\end{equation}
\end{defn}

\begin{prop}[Characterization of $\kappa$-extensible relations]\label{prop:kappa-ext a}
Let $U$ be a relation on $X$. The following conditions are
equivalent.
 \begin{itemize}
\item[(i)] $U$ is $\kappa$-extensible.
\item[(ii)] $cl_{\kappa}(U)$ is a $\kappa$-extension of $U$.
\item[(iii)]$U$ is transitive
and $\nexists$ $x,y,z,t\in X$ such that
\begin{align}
 x\U y\; &\; z\nU y \label{eq:kappa-condition}\\
 x\nU t\;&\; z\U t.\nonumber
 \end{align}
\end{itemize}
\end{prop}

For instance, if we consider the relation $U=\{(x,y),(y,z)\}$ on
$X=\{x,y,z\}$ given above, we have $x\U y$, $y\nU y$, $x\nU z$ and
$y\U z$ and thus, we recover that $U$ has no $\kappa$-extension. One
can also check that the relation "$\mid$\," ("divide") (on $X=[r]$)
defined by $x\mid y$ if and only if "$x$ divide $y$" (i.e.
$\frac{y}{x}\in \mathbb{Z}$) has no $\kappa$-extension. Indeed, the
elements 3,9,2,4 satisfy $3\mid9$,
$2\nmid 9$, $3\nmid 4$ and $2\mid 4$.\\

\pf Clearly (ii)$\implies$(i).

(i)$\implies$(iii): Suppose $U$ has a $\kappa$-extension $S$. Then
\begin{itemize}
\item[(a)]$U$ is transitive: Indeed, let $x,y,z\in X$ and suppose $x\U y$
and $y\U z$. We want to show that $x\U z$. Suppose $x\nU z$, then
since $S$ is a $\kappa$-extension of $U$ and $yUz$, it follows that
$y\S x$ and $x\nS y$. We thus have $x\U y$ and $x\nS y$, which is
impossible since $U\subseteq S$. Thus $x\U z$.
\item[(b)] $\nexists\, x,y,z,t\in X$ satisfying $x\U y$, $z\nU y$, $z\U t$ and $x\nU t$:
Indeed, suppose the contrary.  Then, $x\U y$ and $z\nU y$ imply that
$x\S z$, while $z\U t$ and $x\nU t$ imply that $x\nS z$. We thus
have $x\S y$ and $x\nS y$, which is impossible.
\end{itemize}

(iii)$\implies$(ii): Suppose $U$ satisfy (iii). We want to show that
$H:=cl_{\kappa}(U)$ is a $\kappa$-extension of $U$, that is for any
$x,y,z$ satisfying $x \U y$ and $z\nU y$, we have $x\H z$ and $z\nH
x$. Let $x,y,z$ satisfying $x \U y$ and $z\nU y$. First, by
definition of $H$, we have  $x\H z$. It then  remains to show that
$z\nH x$. Suppose the contrary, i.e. $z\H x$. We distinct two cases:
\begin{itemize}
\item[(a)] $z\U x$: since $U$ is transitive and $x\U y$, we have $z\U y$, which contradicts $z\nU y$.
\item[(b)] $z\nU x$ and $z\H x$: by definition of $H$, there exists $t$
such that $z\U t$ and $x\nU t$. We thus four elements $x,y,z,t$
satisfying $x\U y, z\nU y$, $x\nU t$ and $z\U t$, which contradicts
(iii).
\end{itemize}
 \qed

The following proposition gives some properties of the
$\kappa$-closure.

\begin{prop}\label{prop:kappa-clos}
Let $U$ be a $\kappa$-extensible relation on $X$. Then,
\begin{itemize}
\item $cl_{\kappa}(U)$ is the smallest $\kappa$-extension of $U$ (by
inclusion), i.e. every $\kappa$-extension of $U$ contains
$cl_{\kappa}(U)$.
\item $cl_{\kappa}(U)$ is a bipartitional relation.
\end{itemize}
\end{prop}

\pf  The first assumption is evident by definition of
$cl_{\kappa}(U)$. Set $H:=cl_{\kappa}(U)$. We claim that $H$ is
transitive. Indeed, let $x,y,z\in X$ satisfy $x\H y$ and $y\H z$. We
want to show that $x\H z$. We distinct four cases:
\begin{itemize}
\item[(i)]$x\U y$, $y\U z$: then, by transitivity of $U$, we have
$x\U z$ and thus $x\H z$ (since $U\subseteq H$).
\item[(ii)]$x\U y$, $y\nU z$ and $y\H z$: then by
definition of $H$, there is $t\in X$ such that $y\U t$ and $z\nU t$.
By transitivity of U, we have $x\U t$. We thus have $x\U t$ and
$z\nU t$, which imply, by definition of $H$, that $x\H z$.
\item[(iii)]$x\nU y$ and $x\H y$, $y\U z$: then by definition of $H$, there is $t\in
X$ such that $x\U t$ and $y\nU t$. Suppose $x\nU z$, then the
elements $x,t,y,z$ satisfy $x\U t,\,y\nU t,\,x\nU z,\,y\U z$, which
contradict \eqref{eq:kappa-condition}. We thus have $x\U z$, and in
particular, $x\H z$.
\item[(iv)]$x\nU y$ and $x\H y$, $y\nU z$ and $y\H z$: by definition of $H$,
there exist $t,v\in X$ such that $x\U t$, $y\nU t$, $y\U v$ and
$z\nU v$. Suppose $x\nU v$, then the elements $x,t,y,v$ satisfy
 $x\U t$, $y\nU t$, $x\nU v$ $y\nU v$, which contradict \eqref{eq:kappa-condition}.
Thus we have $x\U v$, and since $z\nU v$, we have by definition of
$H$ that $xHz$.
\end{itemize}
\qed

Let $V$ be a bipartitional relation on $X$ and
$(B_1,\ldots,B_k),(\beta_1,\ldots,\beta_k)$ be the bipartition
associated to $V$ (see Definition~\ref{defn:bipartitional}). Suppose
the block $B_l$ consists of the integer $i_1,i_2,\ldots,i_p$. It
will be convenient to write $c(B_l)$ for the sequence $c(i_1),
c(i_2),\ldots,c(i_p)$ and $m(B_l)=m_l$ for the sum $c(i_1)+
c(i_2)+\ldots+c(i_p)$. In particular, ${m_l\choose c(B_l)}$ will
denote the multinomial coefficient
${c(i_1)+c(i_2)+\ldots+c(i_p)\choose c(i_1),c(i_2),\ldots,c(i_p)}$.

\begin{prop}\label{prop:dist kappa-ext}
Let $U$ be a $\kappa$-extensible relation. It follows that
$H:=cl_{\kappa}(U)$ is a bipartitional relation. Let
$(B_1,\ldots,B_k),(\beta_1,\ldots,\beta_k)$ be the bipartition
associated to $H$. Then,
\begin{align}\label{eq:maj-inv clos}
 \sum_{\w\in\R(\c)}q^{(\maj'_{U}+\inv'_{H\setminus U})(\w)}={c(1)+c(2)+\cdots+c(r)\brack
m_1,m_2,\cdots,m_k}_{q}\;\prod_{l=1}^k{m_l\choose
c(B_l)}q^{\beta_l{m_l\choose c(B_l)}}.
\end{align}
More generally, Equation~\eqref{eq:maj-inv clos} hold for each
 relation $H$ satisfying (1) $H$ is a
$\kappa$-extension of $U$ and (2) $H$ is bipartitional on $X$.
\end{prop}

\pf It is just a combination of Theorem~\ref{thm:majinv},
Proposition~\ref{prop:kappa-clos} and Proposition~2.1
in~\cite{FoZe}.\qed

\section{Proof of Theorem~\ref{thm:maj-inv countable}}

Theorem~\ref{thm:mahonian maj-inv} lead to the following question:
\emph{Given a total order $S$ on $X$, which are the relations $U$ on
$X$ such that $S$ is a $\kappa$-extension of $U$?}

\begin{prop}\label{thm: kappa-extension countable}
Let $U$ be a relation on $X$. The following conditions are
equivalent.
\begin{itemize}
 \item[(i)]The natural order "$>$" is a $\kappa$-extension of $U$.
\item[(ii)] There exists a map $g:X\mapsto X\cup\{\infty\}$
satisfying $g(y)>y$ for each $y\in X$ such that
\begin{align*}
x\U y\Leftrightarrow x\geq g(y).
\end{align*}
\end{itemize}

Moreover, if $U$ satisfy the condition $(i)$, the map $g$ is unique
and defined by
$$g(y)=\left\{
  \begin{array}{ll}
    \min(\{x: x\U y\}), & \hbox{if $\exists x$ such that $x\U y$;} \\
    \infty, & \hbox{otherwise.}
  \end{array}
\right.$$
\end{prop}

\pf (i)$\implies$ (ii): Suppose "$>$" is a $\kappa$-extension of $U$
and let $y\in X$. Then, define $g(y)\in X\cup\{\infty\}$ by
\begin{itemize}
\item $g(y)=\min(\{x ; x\U y\})$ if $\exists x\in X$ satisfying $x\U
y$,
\item $g(y)=\infty $ otherwise.
\end{itemize}

It is clear that $g(y)>y$ for each $y\in X$ because $U \subseteq
">"$.  Let $x,y\in X$. By definition of $g(y)$, we have $x\U
y\implies x\geq g(y)$. Now, suppose $x\geq g(y)$. Since $x\in X$, it
follows that $g(y)<\infty$ and thus, there exists $z\in X$ such that
$z\U y$. We can take $z=g(y)$. Suppose $x\nU y$. Since $z\U y$ and
$x\nU y$ and $">"$ is a $\kappa$-extension of $U$, we then have
$z=g(y)>x$ which contradicts the fact that $x\geq g(y)$. It then
follows that $x\U y$. We thus have proved that
$ x\geq g(y)\implies x\U y$.\\

(ii)$\implies$(i): Let $x,y,z\in X$ satisfying $x\U y$ and $z\nU y$.
It follows from (ii) that $x\geq g(y)$ and $z<g(y)$, and thus $x>z$.
We thus have proved that $">"$ is a $\kappa$-extension of $U$.

\qed

The proof of the following result is left to the reader.
\begin{lem}\label{lem:isomorphime d'ordre}
 The $\kappa$-extensibility on $X$ is
transposable by order isomorphism.

In other words, if $S$ and $T$ are two total orders on $X$ and $h$
is the unique order isomorphism $h:(X,S)\mapsto (X,T)$, i.e $h$ is a
permutation of $X$ and $ x\S y \Leftrightarrow h(x)\T\,h(y)$. Then
$S$ is a $\kappa$-extension of a relation $\U$ on $X$ if and only if
the total order $T$ is a $\kappa$-extension of the relation
$V:="h(U)"$ defined by $x\V y\Leftrightarrow h^{-1}(x)\U h^{-1}(y)$.
\end{lem}

Combining the above Lemma and Proposition~\ref{thm: kappa-extension
countable}, we get immediately the following result.

\begin{prop}\label{prop: preuve de maj-inv countable}
 Let $S$ be a total order on $X$ and $U$ be a relation on $X$.
We denote by $f$ be the (unique) order isomorphism from $(X,">")$ to
$(X,S)$. The following conditions are equivalent.
\begin{itemize}
 \item[(i)]$S$ is a $\kappa$-extension of $U$.
\item[(ii)] There exist an unique map $g:X\mapsto X\cup\{\infty\}$
satisfying $g(y)>y$ for each $y\in X$ such that
\begin{align*}
x\U y\Leftrightarrow f(x)\geq g(f(y)).
\end{align*}
\end{itemize}
\end{prop}

Clearly, Theorem~\ref{thm:maj-inv countable} is an immediate
consequence of Theorem~\ref{thm:mahonian maj-inv} and
Proposition~\ref{prop: preuve de maj-inv countable}.

%%%%****%%%%****%%%%****%%%%****%%%%****%%%%****%%%%****%%%%*
%%%%
%%%% New Section
%%%%
%%%%****%%%%****%%%%****%%%%****%%%%****%%%%****%%%%****%%%%*

\section{Applications: new mahonian statistics}

 In this section, we give some examples of mahonian maj-inv statistics on
$X^*$ which can be derived from the results obtained in this paper.
Such statistics are entirely characterized in
Theorem~\ref{thm:mahonian maj-inv} and Theorem~\ref{thm:maj-inv
countable}.\\

Let $g_k$, $k\in[1,\infty[$, be the maps $X\mapsto X\cup\{\infty\}$
defined for $x\in X$ by
\begin{align*}
g_k(x)&=\lfloor kx+1 \rfloor .\chi(kx<r)+\infty.\chi(kx\geq r)\,.
\end{align*}
Clearly, for each $x\in X$, we have $g_k(x)>x$. By applying
Corollary~\ref{cor:mah maj-inv id} (or Theorem~\ref{thm:maj-inv
countable} with $f=Id$), we obtain immediately the following result.

\begin{prop}\label{prop:mah stat A}

 The statistics $\stat_{g_k}$, $k\in[1,\infty[$, defined for
$\w=x_1x_2\ldots x_n\in X^*$ by
\begin{align*}
{\stat}_{g_k}(\w)&=\sum_{i=1}^{n-1}i.\chi(\frac{x_i}{x_{i+1}}> k)
+\sum_{1\leq i<j\leq n}\chi(k\geq\frac{x_i}{x_{j}}>1)
\end{align*}
are mahonian on $X^*$.
\end{prop}

Note that $\stat_{g_1}=\inv$ and $\stat_{g_r}=\maj$. Now for each
$B\subseteq X$, let $H_B:X\mapsto X\cup\{\infty\}$ be the map
defined for $x\in X$ by $H_B(x)=(x+1).\chi(x\in B,\,x\neq
r)+\infty.\chi(x\notin B \;\text{or}\;x=r )$. Since $H_B(x)>x$ for
each $x\in X$, we obtain by applying Corollary~\ref{cor:mah maj-inv
id} the following result.

\begin{prop}

The statistics $\stat_{H_{B}}$, $B\subseteq X$, defined for
$\w=x_1x_2\ldots x_n\in X^*$ by
\begin{align*}
{\stat}_{H_{B}}(\w)=&\sum_{i=1}^{n-1}i.\chi(x_i>x_{i+1},\;\text{$x_{i+1}\in
B$}) + \sum_{1\leq i<j\leq n}\chi( x_i>x_{j},\;\text{$x_{j}\notin
B$})\end{align*} are mahonian on $X^*$.
\end{prop}

For instance, if $B=\{\text{even numbers $\leq r$} \}$, then the
statistic $\stat_{H_{B}}$ defined for words $\w=x_1x_2\ldots x_n\in
X^*$ by
 \begin{align*}
{\stat}_{H_{B}}(\w)=&\sum_{i=1}^{n-1}i.\chi(x_i>x_{i+1},\;\text{$x_{i+1}$
is even}) + \sum_{1\leq i<j\leq n}\chi( x_i>x_{j},\;\text{$x_{j}$ is
odd})
\end{align*}
is mahonian on $X^*$.

More generally, for $A,B\subseteq X$, Let $U_{A,B}$ be the relation
on $X$ defined by
$$
(x,y)\in U_{A,B}\Longleftrightarrow x\in A\;,\;y\in B\;\text{and}\;
x>y\,.
$$

Suppose $(x,y),(y,z)\in U_{A,B}$. By definition of $U_{A,B}$, we
have $x,y\in A$, $y,z\in B$ and $x>y$ and $y>z$. In particular,
$x\in A$, $z\in B$ and $x>z$, i.e. $(x,z)\in U_{A,B}$. It follows
that $U_{A,B}$ is transitive. Now suppose there exist $x,y,z,t\in X$
such that $(x,y),(z,t)\in U_{A,B}$ and $(x,t),(z,y)\notin U_{A,B}$.
By definition of $U_{A,B}$, we have $x,z\in A$, $y,t\in B$ and
$x>y$, $z\leq y$, $x\leq t$, $z>t$. In particular, $x\leq t$ and
$x>t$, which is impossible. It then follows from
Proposition~\ref{prop:kappa-ext a}  that $U_{A,B}$ is a
$\kappa$-extensible relation on $X$. Let $S_{A,B}$ and $S'_{A,B}$ be
the relation defined on $X$ by
\begin{align*}
S_{A,B}&=\{(x,y)\in X^2\,|\,x\in A,\;y\notin A\}\cup \{(x,y)\in
A^2\,|\,x>y\}\\
 S'_{A,B}&=S_{A,B}\cup \{(x,y)\in (A^c)^2\,|\,x>y\}.
\end{align*}
Then the reader can check that $S'_{A,B}$ and $S_{A,B}$ are two
$\kappa$-extensions of $U_{A,B}$. Suppose
$A=\{a_1,a_2,\ldots,a_k\}_>$. It is then easy to see that $S_{A,B}$
is a bipartitional relation and its associated bipartition is the
pair composed by the partition
$(\{a_1\},\{a_2\},\ldots,\{a_k\},A^c),$ and the null vector
$\textbf{0}=(0,0,\ldots,0,0)$.

Set $\stat_{A,B}:={\maj}'_{U_{A,B}}+{\inv}'_{S_{A,B}\setminus
U_{A,B}}$ and
$\stat'_{A,B}:={\maj}'_{U_{A,B}}+{\inv}'_{S'_{A,B}\setminus
U_{A,B}}$. By definition, the statistics ${\stat}_{A,B}$ and
${\stat}'_{A,B}$ are defined on words $\w=x_1\ldots x_n \in X^*$ by
\begin{align*}
{\stat}_{A,B}(\w)&=\sum_{i=1}^{n-1}i.\chi(x_i>x_{i+1},\,x_i\in
A,\,x_{i+1} \in B)+ \sum_{1\leq i<j\leq n}\chi(x_i>x_j,\,x_{i}\in
A,\,x_j\in A\setminus B) \\
&+\sum_{1\leq i<j\leq n}\chi(x_i\leq x_j,\,x_{i}\in A,\,x_j\in
B\setminus A)+ \sum_{1\leq i<j\leq n}\chi(x_{i}\in A,\,x_j\notin
A\cup
B),\\
{\stat}'_{A,B}(\w)&={\stat}_{A,B}(\w)+ \sum_{1\leq i<j\leq
n}\chi(x_i>x_j,\, \,x_{i}\,\text{and}\,x_j\,\notin A\,).\\
\end{align*}

Applying Theorem~\ref{thm:mahonian maj-inv} and
Proposition~\ref{prop:dist kappa-ext}, we obtain the following
result.
\begin{prop}\label{prop:statAB}
The statistic ${\stat}'_{A,B}$ is mahonian on $X^*$ and
 for each $\c$,
\begin{align*}
 \sum_{\w\in\R(\c)}q^{\stat_{A,B}(\w)}={m(A^c)\choose c(A^c)}
\;{c(1)+c(2)+\cdots+c(r) \brack
c(a_1),c(a_2),\ldots,c(a_k),m(A^c)}_{q}.
\end{align*}
\end{prop}

For instance, if $E=$\{even integers $\leq r$\} and  $O=$\{odd
integers $\leq r$\}, then the statistics ${\stat}_{E,O}$ and
${\stat}'_{E,O}$ are defined for $\w=x_1\ldots x_n\in X^*$ by
\begin{align*}
{\stat}_{E,O}(\w)&=\sum_{i=1}^{n-1}i.\chi(x_i>x_{i+1},\,x_i\;
\text{even},\,x_{i+1} \; \text{odd}) + \sum_{1\leq i<j\leq
n}\chi(x_i>x_j,\,x_{i}\,\text{and}\,x_j\;\text{even}\,)\\
&+ \sum_{1\leq i<j\leq n}\chi(x_i\leq
x_j,\,x_{i}\,\text{even},\,x_j\;\text{odd}\,),\\
{\stat}'_{E,O}(\w)&=\sum_{i=1}^{n-1}i.\chi(x_i>x_{i+1},\,x_i\;
\text{even},\,x_{i+1} \; \text{odd}) + \sum_{1\leq i<j\leq
n}\chi(x_i\leq x_j,\,x_{i}\,\text{even},\,x_j\;\text{odd}\,)\\
+& \sum_{1\leq i<j\leq
n}\chi(x_i>x_j,\,x_{i}\,\text{and}\,x_j\;\text{have the same
parity}\,).
\end{align*}
It then follows from Proposition~\ref{prop:statAB} that the
statistic ${\stat}'_{E,O}$ is mahonian  and the generating function
of ${\stat}_{E,O}$ on each $\R(\c)$ is given by
\begin{align*}
 \sum_{\w\in\R(\c)}q^{\stat_{E,O}(\w)}=&{c(1)+c(3)+\ldots+c(2\left\lfloor\frac{r-1}{2}\right\rfloor +1)
\choose c(1),c(3),\ldots,c(2\left\lfloor\frac{r-1}{2}\right\rfloor +1)}\\
&\times \;{c(1)+c(2)+\cdots+c(r) \brack
c(2),c(4),\ldots,c(2\left\lfloor\frac{r}{2}\right\rfloor),c(1)+c(3)+\ldots+c(2\left\lfloor\frac{r-1}{2}\right\rfloor
+1)}_{q}.
\end{align*}

In particular, if $\mathcal{S}_r$ is the symmetric group of order
$r$, then
\begin{align*}
 \sum_{\sigma\in \,\mathcal{S}_r}q^{\stat_{E,O}(\sigma)}=&\left(\left\lfloor\frac{r+1}{2}\right\rfloor\right)\Large{!}\,\times \frac{[r]_q!}{ [(\left\lfloor\frac{r+1}{2}\right\rfloor)]_{q}!}.
\end{align*}
\vspace{0.5cm}

{\bf Acknowledgments.} I am grateful to Jiang Zeng for helpful
discussions and to an anonymous referee for a careful reading of the
paper.

\end{document}